\documentclass[a4paper]{article}
\usepackage{amssymb}
\usepackage{amsmath}
\usepackage{amsthm}
\usepackage{authblk}
\usepackage[USenglish]{babel}
\usepackage{esint}
\usepackage{graphicx}
\usepackage[colorlinks=true]{hyperref}
\usepackage[latin1]{inputenc}
\newtheorem{theorem}{Theorem}[section]
\newtheorem{lemma}[theorem]{Lemma}
\newtheorem{proposition}[theorem]{Proposition}
\newtheorem{corollary}[theorem]{Corollary}

\newtheorem{definition}[theorem]{Definition}

\newcommand{\N}{\mathbb N}

\newcommand{\R}{\mathbb R}
\newcommand{\C}{\mathbb C}
\renewcommand{\S}{\mathbb S}
\renewcommand{\P}{\mathbb P}

\newcommand{\mcal}{\mathcal}

\newcommand{\mrm}{\mathrm}

\renewcommand{\a}{\alpha}

\renewcommand{\d}{\delta}
\newcommand{\D}{\Delta}
\newcommand{\e}{\varepsilon}
\newcommand{\z}{\zeta}
\renewcommand{\t}{\theta}

\newcommand{\la}{\lambda}
\newcommand{\La}{\Lambda}
\newcommand{\s}{\sigma}

\newcommand{\Si}{\Sigma}

\newcommand{\ph}{\varphi}

\renewcommand{\O}{\Omega}
\newcommand{\wt}{\widetilde}

\newcommand{\ol}{\overline}

\newcommand{\mri}{\mathring}

\newcommand{\fr}{\frac}
\newcommand{\pa}{\partial}
\newcommand{\n}{\nabla}
\newcommand{\fa}{\forall}
\newcommand{\ex}{\exists}
\newcommand{\es}{\emptyset}
\newcommand{\wk}{\rightharpoonup}

\newcommand{\us}{\underset}

\newcommand{\To}{\Rightarrow}

\newcommand{\lto}{\longrightarrow}

\newcommand{\sm}{\setminus}
\renewcommand{\Cup}{\bigcup}

\newcommand{\sub}{\subset}
\newcommand{\Sub}{\Subset}

\newcommand{\nin}{\not\in}

\newcommand{\x}{\times}

\newcommand{\cd}{\cdot}
\newcommand{\ds}{\dots}

\newcommand{\tx}{\text}
\newcommand{\q}{\quad}
\renewcommand{\l}{\left}
\renewcommand{\r}{\right}

\newcommand{\bthm}{\begin{theorem}}
\newcommand{\ethm}{\end{theorem}}
\newcommand{\blem}{\begin{lemma}}
\newcommand{\elem}{\end{lemma}}
\newcommand{\bprop}{\begin{proposition}}
\newcommand{\eprop}{\end{proposition}}
\newcommand{\bcor}{\begin{corollary}}
\newcommand{\ecor}{\end{corollary}}
\newcommand{\bdefi}{\begin{definition}}
\newcommand{\edefi}{\end{definition}}
\newcommand{\bpf}{\begin{proof}}
\newcommand{\epf}{\end{proof}}
\newcommand{\bl}{\begin{array}{l}}
\newcommand{\bll}{\begin{array}{ll}}
\newcommand{\barr}{\begin{array}}
\newcommand{\earr}{\end{array}}
\newcommand{\bite}{\begin{itemize}}
\newcommand{\eite}{\end{itemize}}
\newcommand{\bequ}{\begin{equation}}
\newcommand{\eequ}{\end{equation}}
\newcommand{\beqa}{\begin{eqnarray}}
\newcommand{\eeqa}{\end{eqnarray}}
\newcommand{\beqy}{\begin{eqnarray*}}
\newcommand{\eeqy}{\end{eqnarray*}}

\newcommand{\qm}[1]{``#1''}

\begin{document}

\everymath{\displaystyle}

\title{Moser-Trudinger inequalities for singular Liouville systems}
\author{Luca Battaglia\thanks{Universit\'e Catholique de Louvain, Institut de Recherche en Math\'ematique et Physique, Chemin du Cyclotron $2$, $1348$ Louvain-la-Neuve (Belgium) - luca.battaglia@uclouvain.be\\The author has been supported by the PRIN project \emph{Variational and perturbative aspects of nonlinear differential problems}.$2010$ Mathematics Subject Classification: $35$A$23$, $35$J$47$, $35$J$50$, $58$E$35$.}}
\date{}

\maketitle\

\begin{abstract}
\noindent In this paper we prove a Moser-Trudinger inequality for the Euler-Lagrange functional of general singular Liouville systems on a compact surface. We characterize the values of the parameters which yield coercivity for the functional, hence the existence of energy-minimizing solutions for the system, and we give necessary conditions for boundedness from below.\\
We also provide a sharp inequality under assuming the coefficients of the system to be non-positive outside the diagonal.\\
The proofs use a concentration-compactness alternative, Poho\v zaev-type identities and blow-up analysis.
\end{abstract}\

\section{Introduction}\

An essential tool in the study of the embeddings of Sobolev spaces is the Moser-Trudinger inequality, which gives compact embedding in any $L^p$ space for finite $p\ge1$ and also exponential integrability.\\
If we consider a $2$-dimensional compact Riemannian manifold $(\Si,g)$, due to well-known works from Moser \cite{mos} and Fontana \cite{fon} we get
\bequ\label{mtscal}
\log\int_\Si e^u\mrm dV_g-\fint_\Si u\mrm dV_g\le\fr{1}{16\pi}\int_\Si|\n u|^2\mrm dV_g+C\q\q\q\fa\,u\in H^1(\Si),
\eequ
where $\n=\n_g$ is the gradient given by the metric $g$ and $C=C_{\Si,g}$ is a constant depending only on $\Si$ and $g$.\\

This inequality has fundamental importance in the study of the Liouville equations of the kind
\bequ\label{liour}
-\D u=\rho\l(\fr{he^u}{\int_\Si he^u\mrm dV_g}-1\r),
\eequ
where $\D=\D_g$ is the Laplace-Beltrami operator, $\rho$ a positive real parameter, $h$ a positive smooth function and $\Si$ is supposed, without loss of generality, to have area equal to $|\Si|=1$. In fact, the solutions of \eqref{liour} are critical points of the functional
$$I_\rho(u)=\fr{1}2\int_\Si|\n u|^2\mrm dV_g-\rho\l(\log\int_\Si he^u\mrm dV_g-\int_\Si u\mrm dV_g\r);$$
using the inequality \eqref{mtscal} we can control the last term by the Dirichlet energy, thus showing that $I_\rho$ is bounded from below on $H^1(\Si)$ if and only if $\rho$ is smaller or equal to $8\pi$.\\
Equations like \eqref{liour} have great importance in different contexts like the Gaussian curvature prescription problem (see for instance \cite{cy87,cy88}) and abelian Chern-Simons models in theoretical physics (\cite{tar08,yang}).\\

An extension of the inequality \eqref{mtscal}, which takes into consideration power-type weights, was given by Chen \cite{chenwx} and Trojanov \cite{tro}. For a given $p\in\Si$ and $\a\in(-1,0]$, they showed that
\bequ\label{mtsing}
(1+\a)\l(\log\int_\Si d(\cd,p)^{2\a}e^u\mrm dV_g-\int_\Si u\mrm dV_g\r)\le\fr{1}{16\pi}\int_\Si|\n u|^2\mrm dV_g+C\q\q\q\fa\,u\in H^1(\Si).
\eequ
This inequality allows to treat singularities in the equation \eqref{liour}, that is to consider equations like
\bequ\label{lious}
-\D u=\rho\l(\fr{he^u}{\int_\Si he^u\mrm dV_g}-1\r)-4\pi\sum_{m=1}^M\a_m(\d_{p_m}-1),
\eequ
where we take arbitrary $p_1,\ds,p_M\in\Si$ and $\a_m>-1$ for any $m\in\{1,\ds,M\}$.\\
This is a natural extension of \eqref{liour}, which allows to consider the same problems in a more general context. For instance, it arises in the Gaussian curvature prescription problem on surfaces with conical singularities and in Chern-Simons vortices theory.\\
Defining $G_p$ as the Green function of $-\D$ on $\Si$ centered at a point $p$, through the change of variables
\bequ\label{cov}
u\mapsto u+4\pi\sum_{m=1}^M\a_mG_{p_m}
\eequ
equation \eqref{lious} becomes
$$-\D u=\rho\l(\fr{\wt he^u}{\int_\Si\wt he^u\mrm dV_g}-1\r)$$
with $\wt h=he^{-4\pi\sum_{m=1}^M\a_mG_{p_m}}$.\\
Since $G_p$ has the same behavior as $\fr{1}{2\pi}\log\fr{1}{d(\cd,p)}$ around $p$, then $\wt h$ behaves like $d(\cd,p_m)^{2\a_m}$ around each singular point $p_m$, hence the energy functional 
$$I_\rho(u)=\fr{1}2\int_\Si|\n u|^2\mrm dV_g-\rho\l(\log\int_\Si\wt he^u\mrm dV_g-\int_\Si u\mrm dV_g\r)$$
can be estimated, as in the regular case, using \eqref{mtsing}.\\

The purpose of this paper is to extend inequality \eqref{mtsing} to singular Liouville systems of the type
$$-\D u_i=\sum_{j=1}^Na_{ij}\rho_j\l(\fr{h_je^{u_j}}{\int_\Si h_je^{u_j}\mrm dV_g}-1\r)-4\pi\sum_{m=1}^M\a_{im}(\d_{p_m}-1),\q\q\q i=1,\ds,N,$$
where $A=(a_{ij})$ is a $N\x N$ symmetric positive definite matrix and $\rho_i,h_i,\a_{im}$ are as before.\\
Applying, similarly to \eqref{cov}, the change of variables
$$u_i\mapsto u_i+4\pi\sum_{m=1}^M\a_{im}G_{p_m},$$
the system becomes
\bequ\label{liousy}
-\D u_i=\sum_{j=1}^Na_{ij}\rho_j\l(\fr{\wt h_je^{u_j}}{\int_\Si\wt h_je^{u_j}\mrm dV_g}-1\r),\q\q\q i=1,\ds,N,
\eequ
with $\wt h_j$ having the same behavior around the singular points.\\
The system has a variational formulation with the energy functional
\bequ\label{jrho}
J_\rho(u):=\fr{1}2\sum_{i,j=1}^Na^{ij}\int_\Si\n u_i\cd\n u_j\mrm dV_g-\sum_{i=1}^N\rho_i\l(\log\int_\Si\wt h_ie^{u_i}\mrm dV_g-\int_\Si u_i\mrm dV_g\r),
\eequ
with $a^{ij}$ indicating the entries of the inverse matrix $A^{-1}$ of $A$.\\

A recent paper by the author and Malchiodi (\cite{batmal}) gives an answer for the particular case of the $SU(3)$ Toda system, that is $N=2$ and $A$ is the Cartan matrix
$$\l(\barr{cc}2&-1\\-1&2\earr\r).$$
This is a particularly interesting case, due to its application in the description of holomorphic curves in $\C\P^N$ in geometry (\cite{bw,cal,cw}) and in the non-abelian Chern-Simons theory in physics (\cite{dunne,tar08,yang}).\\
The authors prove a sharp inequality, that is they show that the functional $J_\rho$ is bounded from below if and only if both the parameters $\rho_i$ are less or equal than $4\pi\min\l\{1,1+\min_m\a_{im}\r\}$, thus extending the result in the regular case from \cite{jw01}.\\

Concerning general regular Liouville systems, Wang \cite{wang} gave necessary and sufficient conditions for the boundedness from below of $J_\rho$, following previous results in \cite{csw97,csw02} for the problem on Euclidean domains with Dirichlet boundary conditions. Analogous results were given in \cite{sw1} for the standard unit sphere $\l(\S^2,g_0\r)$ and in \cite{sw2} for a similar problem.\\
In these papers, the authors introduce, for any $\mcal I\sub\{1,\ds,N\}$, the following function of the parameter $\rho$:
$$\La_{\mcal I}(\rho)=8\pi\sum_{i\in\mcal I}\rho_i-\sum_{i,j\in\mcal I}a_{ij}\rho_i\rho_j.$$
What they prove is boundedness from below for $J_\rho$ for any $\rho\in\R_+^N$ which satisfies $\La_{\mcal I}(\rho)>0$ for all the subsets $\mcal I$ of $\{1,\ds,N\}$, whereas $\inf_{H^1(\Si)^N}J_\rho=-\infty$ whenever $\La_{\mcal I}(\rho)<0$ for some $\mcal I\sub\{1,\ds,N\}$.\\

The first main result of this paper is an extension of the results from \cite{csw97,csw02,wang} to the case of singularities.\\
Similarly to Liouville equation, we will have to multiply some quantities by $1+\a_{im}$. Precisely, we have:\\

\bthm\label{subcrit}${}$\\
Let $J_\rho$ be the functional defined by \eqref{jrho} and set, for $\rho\in\R_{>0}^N,\,x\in\Si$ and $i\in\mcal I\sub\{1,\ds,N\}$:
\bequ\label{alfai}
\a_i(x)=\l\{\bll\a_{im}&\tx{if }x=p_m\\0&\tx{otherwise}\earr\r.\q\q\q\La_{\mcal I,x}(\rho):=8\pi\sum_{i\in\mcal I}(1+\a_i(x))\rho_i-\sum_{i,j\in\mcal I}a_{ij}\rho_i\rho_j
\eequ
$$\La(\rho):=\min_{\mcal I\sub\{1,\ds,N\},x\in\Si}\La_{\mcal I,x}(\rho).$$
Then, $J_\rho$ is bounded from below if $\La(\rho)>0$, whereas $J_\rho$ is unbounded from below if $\La(\rho)<0$.
\ethm\

Notice that, in the definition of $\La$, the minimum makes sense because it is taken in a finite set, since $\a_i(x)=0$ for all points of $\Si$ but a finite number, and for all the former points $\La_{\mcal I,x}$ is defined in the same way.\\
As a consequence of this result, we obtain information about the existence of solutions for the system \eqref{liousy}.\\

\bcor\label{ex}${}$\\
The functional $J_\rho$ is coercive in $\ol H^1(\Si)$ if and only if $\La(\rho)>0$.\\
Therefore, if this occurs, then $J_\rho$ admits a minimizer $u$ which solves \eqref{liousy}.
\ecor\

Theorem \ref{subcrit} leaves an open question about what happens when $\La(\rho)=0$. In this case, as we will see in the following Sections, one encounters blow-up phenomena which are not yet fully known for general systems.\\
Anyway, we can say something more if we assume in addition $a_{ij}\le0$ for any $i,j\in\{1,\ds,N\}$ with $i\ne j$. First of all, we notice that in this case $$\La(\rho)=\min_{i\in\{1,\ds,N\}}\l(8\pi(1+\wt\a_i)\rho_i-a_{ii}\rho_i^2\r),\q\q\q\tx{where}$$
\bequ\label{alfatilde}
\wt\a_i:=\min_{m\in\{1,\ds,M\},x\in\Si}\a_i(x)=\min\l\{0,\min_{m\in\{1,\ds,M\}}\a_{im}\r\};
\eequ
hence the sufficient condition in Theorem \ref{subcrit} is equivalent to assuming $\rho_i<\fr{8\pi(1+\wt\a_i)}{a_{ii}}$ for any $i$.\\
With this assumption, studying what happens when $\La_{\mcal I}(\rho)=0$ is reduced to a single-component local blow-up, which can be treated by using an inequality from \cite{as}. Therefore, we get the following sharp result:\\

\bthm\label{sharp}${}$\\
Let $J_\rho$ be defined by \eqref{jrho}, $\wt\a_i$ as in \eqref{alfatilde} and $\La(\rho)$ as in Theorem \ref{subcrit}, and suppose $a_{ij}\le0$ for any $i,j\in\{1,\ds,N\}$ with $i\ne j$.\\
Then, $J_\rho$ is bounded from below on $H^1(\Si)^N$ if and only if $\La(\rho)\ge0$, namely if and only if $\rho_i\le\fr{8\pi(1+\wt\a_i)}{a_{ii}}$ for any $i\in\{1,\ds,N\}$.
\ethm\

We remark that the assuming $A$ to be positive definite is necessary. If it is not, then $J_\rho$ is unbounded from below for any $\rho$.\\
In fact, suppose there exists $v\in\R^N$ such that $\sum_{i,j=1}^Na^{ij}v_iv_j\le-\t|v|^2$ for some $\t>0$. Then, we consider the family of functions $u^\la(x):=\la v\cd x$; by Jensen's inequality we get
\beqy
J_\rho\l(u^\la\r)&\le&\fr{1}2\sum_{i,j=1}^Na^{ij}\int_\Si\n u_i^\la\cd\n u_j^\la\mrm dV_g-\sum_{i=1}^N\rho_i\int_\Si\log\wt h_i\mrm dV_g\\
&\le&-\fr{\t}2\la^2|v|^2+C\\
&\us{n\to+\infty}\lto&-\infty.
\eeqy
We also notice that, with respect to the scalar case, in Theorem \ref{subcrit} and Corollary \ref{ex} the positive coefficients $\a_{im}$'s may affect the definition of $\La(\rho)$, hence the conditions for coercivity and boundedness from below of $J_\rho$.\\
On the other hand, under the assumptions of Theorem \ref{sharp}, coercivity and boundedness from below only depend on the negative $\a_{im}$'s, just like for the scalar equation.\\

The plan of this paper is the following: in Section $2$ we will introduce some notations and some preliminary results which will be used throughout the rest of the paper. In Section $3$ we will show a sort of Concentration-compactness theorem, showing the possible non-compactness phenomena for solutions of the system \eqref{liousy}. Finally, in Sections $4$ and $5$ we will give the proof of the two main theorems.\\

\section{Notations and preliminaries}\

In this section, we will give some useful notation and some known preliminary results which will be needed to prove the two main theorems.\\

Given two points $x,y\in\Si$, we will indicate the metric distance on $\Si$ between them as $d(x,y)$. We will indicate the open metric ball centered in $p$ having radius $r$ as
$$B_r(x):=\{y\in\Si:\,d(x,y)<r\}.$$
For any subset of a topological space $A\sub X$ we indicate its closure as $\ol A$ and its interior part as $\mri A$.\\

Given a function $u\in L^1(\Si)$, the symbol $\ol u$ will indicate the average of $u$ on $\Si$. Since we assume $|\Si|=1$, we can write:
$$\ol u=\int_\Si u\mrm dV_g=\fint_\Si u\mrm dV_g.$$
We will indicate the subset of $H^1(\Si)$ which contains the functions with zero average as
$$\ol H^1(\Si):=\l\{u\in H^1(\Si):\,\ol u=0\r\}.$$
Since the functional $J_\rho$ defined by \eqref{jrho} is invariant by addition of constants, it will not be restrictive to study it on $\ol H^1(\Si)^N$ rather than on $H^1(\Si)^N$.\\

We will indicate with the letter $C$ large constants which can vary among different lines and formulas. To underline the dependence of $C$ on some parameter $\a$, we indicate with $C_\a$ and so on.\\
We will denote as $o_\a(1)$ quantities which tend to $0$ as $\a$ tends to $0$ or to $+\infty$ and we will similarly indicate bounded quantities as $O_\a(1)$, omitting in both cases the subscript(s) when it is evident from the context.\\

First of all, we need a result from Brezis and Merle \cite{bremer}. It is a classical estimate about exponential integrability of solutions of some elliptic PDEs.\\

\blem(\cite{bremer}, Theorem $1$)\label{bremer1}\\
Take $r>0,\,\O:=B_r(0)\sub\R^2,\,f\in L^1(\O)$ with $\|f\|_{L^1(\O)}<4\pi$ and $u$ solving
$$\l\{\bll-\D u=f&\tx{in }\O\\u=0&\tx{on }\pa\O\earr\r..$$
Then, for any $q\in\l[1,\fr{4\pi}{\|f\|_{L^1(\O)}}\r)$ there exists a constant $C=C_{q,r}$ such that $\int_\O e^{q|u(x)|}\mrm dx\le C$.
\elem\

A crucial role in the proof of both Theorem \ref{subcrit} and \ref{sharp} will be played by the concentration values of the sequences of solutions of \eqref{liousy}.\\
For a sequence $u^n=\l\{u_1^n,\ds,u_N^n\r\}_{n\in\N}$ of solutions of \eqref{liousy} with $\rho=\rho^n=\l\{\rho_1^n,\ds,\rho_N^n\r\}$, we define (up to subsequences), for $i\in\{1,\ds,N\}$, the concentration value of its $i^{th}$ component around a point $x\in\Si$ as
\bequ\label{sigmai}
\s_i(x):=\lim_{r\to0}\lim_{n\to+\infty}\rho_i^n\fr{\int_{B_r(x)}\wt h_ie^{u_i^n}\mrm dV_g}{\int_\Si\wt h_ie^{u_i^n}\mrm dV_g}.
\eequ
In a recent paper (\cite{lwzhang}, see also \cite{jlw} for the regular case) it was proved, by a Poho\v zaev identity, that the concentration values satisfy the following algebraic relation, which involves the same quantities as in Theorem \ref{subcrit}:\\

\bprop(\cite{jlw}, Lemma $2.2$; \cite{lwzhang}, Proposition $3.1$)\label{lwzhang}\\
Let $\l\{u^n\r\}_{n\in\N}$ be a sequence of solutions of \eqref{liousy}, $\a_i(x)$ and $\La_{\mcal I,x}$ as in \eqref{alfai} and $\s(x)=(\s_1(x),\ds,\s_N(x))$ as in \eqref{sigmai}. Then,
$$\La_{\{1,\ds,N\},x}(\s(x))=8\pi\sum_{i=1}^N(1+\a_i(x))\s_i(x)-\sum_{i,j=1}^Na_{ij}\s_i(x)\s_j(x)=0.$$
\eprop\

To study the concentration phenomena of solutions of \eqref{liousy} we will use the following simple but useful calculus Lemma:\\

\blem(\cite{jw01}, Lemma $4.4$)\label{fanbn}\\
Let $\l\{a^n\r\}_{n\in\N}$ and $\l\{b^n\r\}_{n\in\N}$ two sequences of real numbers satisfying
$$a^n\us{n\to+\infty}\lto+\infty\q\q\q\lim_{n\to+\infty}\fr{b^n}{a^n}\le0.$$
Then, there exists a smooth function $F:[0,+\infty)\to\R$ which satisfies, up to subsequences,
$$0<F'(t)<1\q\fa\,t>0\q\q\q F'(t)\us{t\to+\infty}\lto0\q\q\q F\l(a^n\r)-b^n\us{n\to+\infty}\lto+\infty.$$
\elem\

Finally, as anticipated in the introduction, we will need a singular Moser-Trudinger inequality for Euclidean domains by Adimurthi and Sandeep \cite{as}, and its straightforward corollary.\\

\bthm(\cite{as}, Theorem $2.1$)\\
For any $r>0,\,\a\in(-1,0]$ there exists a constant $C=C_{\a,r}$ such that if $\O:=B_r(0)\sub\R^2$ and $u\in H^1_0(\O)$, then
$$\int_\O|\n u(x)|^2\mrm dx\le1\q\To\q\int_\O|x|^{2\a}e^{4\pi(1+\a)u(x)^2}\mrm dx\le C$$
\ethm\

\bcor\label{as}${}$\\
For any $r>0,\,\a\in(-1,0]$ there exists a constant $C=C_{\a,r}$ such that if $\O:=B_r(0)\sub\R^2$ and $u\in H^1_0(\O)$, then
$$(1+\a)\log\int_\O|x|^{2\a}e^{u(x)}\mrm dx\le\fr{1}{16\pi}\int_\O|\n u(x)|^2\mrm dx+C$$
\ecor\

\bpf${}$\\
By the elementary inequality $u\le\t u^2+\fr{1}{4\t}$ with $\t=\fr{4\pi(1+\a)}{\int_\O|\n u(y)|^2\mrm dy}$ we get
\beqy
(1+\a)\log\int_\O|x|^{2\a}e^{u(x)}\mrm dx&\le&(1+\a)\log\int_\O|x|^{2\a}e^{\t u(x)^2+\fr{1}{4\t}}\mrm dx\\
&=&\fr{1}{16\pi}\int_\O|\n u(y)|^2\mrm dy+(1+\a)\log\int_\O|x|^{2\a}e^{4\pi(1+\a)\l(\fr{u(x)}{\sqrt{\int_\O|\n u(y)|^2\mrm dy}}\r)^2}\mrm dx\\
&\le&\fr{1}{16\pi}\int_\O|\n u(y)|^2\mrm dy+C.
\eeqy
\epf\

\section{A Concentration-compactness theorem}\

The aim of this section is to prove a result which describes the concentration phenomena for the solutions of \eqref{liousy}, extending what was done for the two-dimensional Toda system in \cite{batmal,ln}.\\

We actually have to normalize such solutions to bypass the issue of the invariance by translation by constants and to have the parameter $\rho$ multiplying only the constant term.\\
In fact, for any solution $u$ of \eqref{liousy} the functions
\bequ\label{vi}
v_i:=u_i-\log\int_\Si\wt h_ie^{u_i}\mrm dV_g+\log\rho_i
\eequ
solve
\bequ\label{eqvi}
\l\{\bl-\D v_i=\sum_{j=1}^Na_{ij}\l(\wt h_je^{v_j}-\rho_j\r)\\\int_\Si\wt h_ie^{v_i}\mrm dV_g=\rho_i\earr\r.\q\q\q i=1,\ds,N.
\eequ
Moreover, we can rewrite in a shorter way \eqref{sigmai} as
$$\s_i(x)=\lim_{r\to0}\lim_{n\to+\infty}\int_{B_r(x)}\wt h_i^ne^{v_i^n}\mrm dV_g.$$
For such functions, we get the following concentration-compactness alternative:\\

\bthm\label{conc}${}$\\
Let $\{u^n\}_{n\in\N}$ be a sequence of solutions of \eqref{liousy} with $\rho^n\us{n\to+\infty}\lto\rho\in\R_+^N$ and $\wt h_i^n=V_i^n\wt h_i$ with $V_i^n\us{n\to+\infty}\lto1$ in $C^1(\Si)^N$, $\{v^n\}_{n\in\N}$ be defined as in \eqref{vi} and $\mcal S_i$ be defined, for $i\in\{1,\ds,N\}$, by
\bequ\label{si}
\mcal S_i:=\l\{x\in\Si:\,\ex\,x^n\us{n\to+\infty}\lto x\tx{ such that }v_i^n\l(x^n\r)\us{n\to+\infty}\lto+\infty\r\}.
\eequ
Then, up to subsequences, one of the following occurs:
\bite
\item If $\mcal S_i=\es$ for any $i\in\{1,\ds,N\}$, then $v^n\us{n\to+\infty}\lto v$ in $W^{2,q}(\Si)^N$ for some $q>1$ and some $v$ which solves \eqref{eqvi}.
\item If $\mcal S_i\ne\es$ for some $i$, then it is a finite set for all such $i$'s. If this occurs, then there is a subset $\mcal I\sub\{1,\ds,N\}$ such that $v_j^n\us{n\to+\infty}\lto-\infty$ in $L^\infty_{\mrm{loc}}\l(\Si\sm\Cup_{j'=1}^N\mcal S_{j'}\r)$ for any $j\in\mcal I$ and $v_j^n\us{n\to+\infty}\lto v_j$ in $W^{2,q}_{\mrm{loc}}\l(\Si\sm\Cup_{j'=1}^N\mcal S_{j'}\r)$ for some $q>1$ and some suitable $v_j$, for any $j\in\{1,\ds,N\}\sm\mcal I$.\\
\eite
\ethm\

Since $\wt h_j$ is smooth outside the points $p_m$'s, the estimates in $W^{2,q}(\Si)$ are actually in $C^{2,\a}_{\mrm{loc}}\l(\Si\sm\Cup_{m=1}^Mp_m\r)$ and the estimates in $W^{2,q}_{\mrm{loc}}\l(\Si\sm\Cup_{j'=1}^N\mcal S_{j'}\r)$ are actually in $C^{2,\a}_{\mrm{loc}}\l(\Si\sm\l(\Cup_{j'=1}^N\mcal S_{j'}\cup\Cup_{m=1}^Mp_m\r)\r)$. Anyway, estimates in $W^{2,q}$ will suffice in most of the paper.\\

To prove Theorem \ref{conc} we need two preliminary lemmas.\\
The first is a Harnack-type alternative for sequences of solutions of PDEs. It is inspired by \cite{bremer,ln}.\\

\blem\label{har}${}$\\
Let $\O\sub\Si$ be a connected open subset, $\{f^n\}_{n\in\N}$ a bounded sequence in $L^q_{\mrm{loc}}(\O)\cap L^1(\O)$ for some $q>1$ and $\{w^n\}_{n\in\N}$ bounded from above and solving $-\D w^n=f^n$ in $\O$.\\
Then, up to subsequences, one of the following alternatives holds:
\bite
\item $w^n$ is uniformly bounded in $L^\infty_{\mrm{loc}}(\O)$.
\item $w^n\us{n\to+\infty}\lto-\infty$ in $L^\infty_{\mrm{loc}}(\O)$.
\eite
\elem\

\bpf${}$\\
Take a compact set $\mcal K\Sub\O$ and cover it with balls of radius $\fr{r}2$, with $r$ smaller than the injectivity radius of $\Si$. By compactness, we can write $\mcal K\sub\Cup_{h=1}^HB_\fr{r}2(x_h)$. If the second alternative does not occur, then up to relabeling we get $\sup_{B_r(x_1)}w^n\ge-C$.\\
Then, we consider the solution $z^n$ of
$$\l\{\bll-\D z^n=f^n&\tx{in }B_r(x_1)\\z^n=0&\tx{on }\pa B_r(x_1)\earr\r.,$$
which is bounded in $L^\infty(B_r(x_1))$ by elliptic estimates. This means that, for a large constant $C$, the function $C-w^n+z^n$ is positive, harmonic and bounded from below on $B_r(x_1)$, and moreover its infimum is bounded from above; therefore, applying the Harnack inequality (which is allowed since $r$ is small enough) we get that $C-w^n+z^n$ is uniformly bounded in $L^\infty\l(B_\fr{r}2(x_1)\r)$, hence $w^n$ is.\\
At this point, by connectedness, we can relabel the index $h$ in such a way that $B_\fr{r}2(x_h)\cap B_\fr{r}2(x_{h+1})\ne\es$ for any $h\in\{1,\ds,H-1\}$ and we repeat the argument for $B_\fr{r}2(x_2)$: since it has nonempty intersection with $B_\fr{r}2(x_1)$, we have $\sup_{B_r(x_2)}w^n\ge-C$, hence we get boundedness in $L^\infty\l(B_\fr{r}2(x_2)\r)$. In the same way, we obtain the same result in all the balls $B_\fr{r}2(x_h)$, whose union contains $\mcal K$, therefore $w^n$ must be uniformly bounded on $\mcal K$ and we get the conclusion.
\epf\

The second Lemma basically says that if all the concentration values in a point are under a certain threshold, and in particular if all of them equal zero, then compactness occurs around that point.\\
On the other hand, if a point belongs to some set $\mcal S_i$, then at least a fixed amount of mass has to accumulate around it; hence, being the total mass uniformly bounded from above, this can occur only for a finite number of points, so we deduce the finiteness of the $\mcal S_i$'s.\\
Precisely, we have the following, inspired again by \cite{ln}, Lemma $4.4$:\\

\blem\label{linfty}${}$\\
Let $\l\{v^n\r\}_{n\in\N}$ and $\mcal S_i$ be as in \eqref{si} and $\s_i$ as in \eqref{sigmai}, and suppose $\s_i(x)<\s_i^0$ for any $i\in\{1,\ds,N\}$, where
$$\s_i^0:=\fr{4\pi\min\l\{1,1+\min_{j\in\{1,\ds,N\},m\in\{1,\ds,M\}}\a_{jm}\r\}}{\sum_{j=1}^Na_{ij}^+}.$$
Then, $x\nin\mcal S_i$ for any $i\in\{1,\ds,N\}$.
\elem\

\bpf${}$\\
First of all we notice that $\s_i^0$ is well-defined for any $i$ because $a_{ii}>0$, hence $\sum_{j=1}^Na_{ij}^+>0$.\\
Under the hypotheses of the Lemma, for large $n$ and small $r$ we have
\bequ\label{strict}
\int_{B_r(x)}\wt h_i^ne^{v_i^n}\mrm dV_g<\s_i^0.
\eequ
Let us consider $w_i^n$ and $z_i^n$ defined by
\bequ\label{zn}
\l\{\bll-\D w_i^n=-\sum_{j=1}^Na_{ij}\rho_j^n&\tx{in }B_r(x)\\w_i^n=0&\tx{on }\pa B_r(x)\earr\r.,\q\q\l\{\bll-\D z_i^n=\sum_{j=1}^Na_{ij}^+\wt h_j^ne^{v_j^n}&\tx{in }B_r(x)\\z_i^n=0&\tx{on }\pa B_r(x)\earr\r..
\eequ
Is it evident that the $w_i^n$'s are uniformly bounded in $L^\infty(B_r(x))$.\\
As for the $z_i^n$'s, we can suppose to be working on a Euclidean disc, up to applying a perturbation to $\wt h_i^n$ which is smaller as $r$ is smaller, hence for $r$ small enough we still have the strict estimate \eqref{strict}.\\
Therefore, we get
$$\l\|-\D z_i^n\r\|_{L^1(B_r(x))}=\sum_{j=1}^Na_{ij}^+\int_{B_r(x)}\wt h_j^ne^{v_j^n}\mrm dV_g<\sum_{j=1}^Na_{ij}^+\s_j^0\le4\pi\min\{1,1+\a_i(x)\}$$ and we can apply Lemma \ref{bremer1} to obtain $\int_{B_r(x)}e^{q|z_i^n|}\mrm dV_g\le C$ for some $q>\fr{1}{\min\{1,1+\a_i(x)\}}$.\\
If $\a_i(x)\ge0$, then taking $q\in\l(1,\fr{4\pi}{\l\|-\D z_i^n\r\|_{L^1(B_r(x))}}\r)$ we have
$$\int_{B_r(x)}\l(\wt h_i^ne^{z_i^n}\r)^q\mrm dV_g\le C_r\int_{B_r(x)}e^{q|z_i^n|}\mrm dV_g\le C.$$
On the other hand, if $\a_i(x)<0$, we choose
$$q\in\l(1,\fr{4\pi}{\l\|-\D z_i^n\r\|_{L^1(B_r(x))}-4\pi\a_i(x)}\r)\q\q\q q'\in\l(\fr{4\pi}{4\pi-q\l\|-\D z_i^n\r\|_{L^1(B_r(x))}},\fr{1}{-\a_i(x)q}\r)$$ and, applying H\"older's inequality,
\beqy\int_{B_r(x)}\l(\wt h_i^ne^{z_i^n}\r)^q\mrm dV_g&\le&C_r\int_{B_r(x)}d(\cd,x)^{2q\a_i(x)}e^{qz_i^n}\mrm dV_g\\
&\le&C\l(\int_{B_r(x)}d(\cd,x)^{2qq'\a_i(x)}\mrm dV_g\r)^\fr{1}{q'}\l(\int_{B_r(x)}e^{q\fr{q'}{q'-1}|z_i^n|}\mrm dV_g\r)^{1-\fr{1}{q'}}\\
&\le&C,
\eeqy
because $qq'\a_i(x)>-1$ and $q\fr{q'}{q'-1}\a_i(x)<\fr{4\pi}{\l\|-\D z_i^n\r\|_{L^1(B_r(x))}}$. Hence $\wt h_i^ne^{z_i^n}$ is uniformly bounded in $L^q(B_r(x))$ for some $q>1$.\\

Now, let us consider $v_i^n-z_i^n-w_i^n$: it is a subharmonic sequence by construction, so for any $y\in B_\fr{r}2(x)$ we get
\beqy v_i^n(y)-z_i^n(y)-w_i^n(y)&\le&\fint_{B_\fr{r}2(y)}\l(v_i^n-z_i^n-w_i^n\r)\mrm dV_g\\
&\le&C\int_{B_\fr{r}2(y)}(v_i^n-z_i^n-w_i^n)^+\mrm dV_g\\
&\le&C\int_{B_r(x)}\l((v_i^n-z_i^n)^++(w_i^n)^-\r)\mrm dV_g\\
&\le&C\l(1+\int_{B_r(x)}\l(v_i^n-z_i^n\r)^+\mrm dV_g\r).
\eeqy
Moreover, since the maximum principle yields $z_i^n\ge0$, taking $\t=\l\{\bll1&\tx{if }\a_i(x)\le0\\\in\l(0,\fr{1}{1+\a_i(x)}\r)&\tx{if }\a_i(x)>0\earr\r.$, we get
\beqy
\int_{B_r(x)}\l(v_i^n-z_i^n\r)^+\mrm dV_g&\le&\int_{B_r(x)}(v_i^n)^+\mrm dV_g\\
&\le&\fr{1}{e\t}\int_{B_r(x)}e^{\t v_i^n}\mrm dV_g\\
&\le&C\l\|\l(\wt h_i^n\r)^{-\t}\r\|_{L^\fr{1}{1-\t}(B_r(x))}\l(\int_{B_r(x)}\wt h_i^n e^{v_i^n}\mrm dV_g\r)^\t\\
&\le&C.
\eeqy
Therefore, we showed that $v_i^n-z_i^n-w_i^n$ is bounded from above in $B_{\fr{r}2}(x)$, that is $e^{v_i^n-z_i^n-w_i^n}$ is uniformly bounded in $L^\infty\l(B_\fr{r}2(x)\r)$. Since the same holds for $e^{w_i^n}$ and $\wt h_i^ne^{z_i^n}$ is uniformly bounded in $L^q\l(B_\fr{r}2(x)\r)$ for some $q>1$, we deduce that also
$$\wt h_i^ne^{v_i^n}=\wt h_i^ne^{z_i^n}\,e^{v_i^n-z_i^n-w_i^n}\,e^{w_i^n}$$
is bounded in the same $L^q\l(B_\fr{r}2(x)\r)$.\\
Thus, we have an estimate on $\l\|-\D z_i^n\r\|_{L^q\l(B_\fr{r}2(x)\r)}$ for any $i\in\{1,\ds,N\}$, hence by standard elliptic estimates we deduce that $z_i^n$ is uniformly bounded in $L^\infty\l(B_\fr{r}2(x)\r)$. Therefore, we also deduce that
$$v_i^n=\l(v_i^n-z_i^n-w_i^n\r)+z_i^n+w_i^n$$
is bounded from above on $B_\fr{r}2(x)$, which is equivalent to saying $x\nin\Cup_{i=1}^N\mcal S_i$.
\epf\

From this proof, we notice that, under the assumptions of Theorem \ref{sharp}, the same result holds for any single index $i\in\{1,\ds,N\}$. In other words, the upper bound on one $\s_i$ implies that $x\nin\mcal S_i$.\\

\bcor\label{equiv}${}$\\
Suppose $a_{ij}\le0$ for any $i\ne j$.\\
Then, for any given $i\in\{1,\ds,N\}$ the following conditions are equivalent:
\bite
\item $x\in\mcal S_i$.
\item $\s_i(x)\ne0$.
\item $\s_i(x)\ge\s_i'=\fr{4\pi\min\l\{1,1+\min_m\a_{im}\r\}}{a_{ii}}$.
\eite
\ecor\

\bpf${}$\\
The third statement trivially implies the second and the second implies the first, since if $v_i^n$ is bounded from above in $B_r(x)$ then $\wt h_i^ne^{v_i^n}$ is bounded in $L^q(B_r(x))$. Finally, if $\s_i(x)<\s_i'$ then the sequence $\wt h_i^ne^{z_i^n}$ defined by \eqref{zn} is bounded in $L^q$ for $q>1$,so one can argue as in Lemma \ref{linfty} to get boundedness from above of $v_i^n$ around $x$, that is $x\nin\mcal S_i$.
\epf\

We can now prove the main theorem of this Section.\\

\bpf[Proof of Theorem \ref{conc}]${}$\\
If $\mcal S_i=\es$ for any $i$, then $e^{v_i^n}$ is bounded in $L^\infty(\Si)$, so $-\D v_i^n$ is bounded in $L^q(\Si)$ for any $$q\in\l[1,\fr{1}{-\min_{j\in\{1,\ds,N\},m\in\{1,\ds,M\}}\a_{jm}}\r).$$
Therefore, we can apply Lemma \ref{har} to $v_i^n$ on $\Si$, where we must have the first alternative for every $i$, since otherwise the dominated convergence would give $\int_\Si\wt h_i^ne^{v_i^n}\mrm dV_g\us{n\to+\infty}\lto0$ which is absurd; standard elliptic estimates allow to conclude compactness in $W^{2,q}(\Si)$.\\
Suppose now $\mcal S_i\ne\es$ for some $i$; from Lemma \ref{linfty} we deduce
$$|\mcal S_i|\s_i^0\le\sum_{x\in\mcal S_i}\max_j\s_j(x)\le\sum_{j=1}^N\sum_{x\in\mcal S_i}\s_j(x)\le\sum_{j=1}^N\rho_j,$$
hence $\mcal S_i$ is finite.\\
For any $j\in\{1,\ds,N\}$, we can apply Lemma \ref{har} on $\Si\sm\Cup_{j'=1}^N\mcal S_{j'}$ with $f^n=\sum_{j'=1}^Na_{jj'}\l(\wt h_{j'}^ne^{v_{j'}^n}-\rho_{j'}^n\r)$, since the last function is bounded in $L^q_{\mrm{loc}}\l(\Si\sm\Cup_{j'=1}^N\mcal S_{j'}\r)$.\\
Therefore, either $v_j^n$ goes to $-\infty$ or it is bounded in $L^\infty_{\mrm{loc}}$, and in the last case we get compactness in $W^{2,q}_{\mrm{loc}}$ by applying again standard elliptic regularity.\\
\epf\

\section{Proof of Theorem \ref{subcrit}.}\

Here we will prove the theorem which gives sufficient and necessary conditions for the functional $J_\rho$ to be bounded from below.\\
In other words, setting
\bequ\label{e}
E:=\l\{\rho\in\R_+^N:\,J_\rho\tx{ is bounded from below on }H^1(\Si)^N\r\},
\eequ
we will prove that $\l\{\La>0\r\}\sub E\sub\l\{\La\ge0\r\}$.\\

As a first thing, we notice that the set $E$ is not empty and it verifies a simple monotonicity condition.\\

\blem\label{ee}${}$\\
The set $E$ defined by \eqref{e} is nonempty.\\
Moreover, for any $\rho\in E$ then $\rho'\in E$ provided $\rho'_i\le\rho_i$ for any $i\in\{1,\ds,N\}$.
\elem\

\bpf${}$\\
Let $\t>0$ be the biggest eigenvalue of the matrix $(a_{ij})$. Then,
$$J_\rho(u)\ge\sum_{i=1}^N\l(\fr{1}{2\t}\int_\Si|\n u_i|^2\mrm dV_g-\rho_i\l(\log\int_\Si\wt h_ie^{u_i}\mrm dV_g-\ol{u_i}\r)\r).$$
Therefore, from scalar Moser-Trudinger inequality \eqref{mtsing}, we deduce that $J_\rho$ is bounded from below if $\rho_i\le\fr{8\pi(1+\wt\a_i)}\t$, hence $E\ne\es$.\\
Suppose now $\rho\in E$ and $\rho'_i\le\rho_i$ for any $i$. Then, through Jensen's inequality, we get
\beqy
J_{\rho'}(u)&=&J_\rho(u)+\sum_{i=1}^N(\rho_i-\rho'_i)\log\int_\Si e^{u_i-\ol{u_i}+\log\wt h_i}\mrm dV_g\\
&\ge&-C+\sum_{i=1}^N(\rho_i-\rho'_i)\int_\Si\log\wt h_i\mrm dV_g\\
&\ge&-C
\eeqy
for any $u\in H^1(\Si)^N$, hence the claim.
\epf\

It is interesting to observe that a similar monotonicity condition is also satisfied by the set $\{\La>0\}$ (although one can easily see that it is not true if we replace $\La$ with $\La_{\mcal I,x}$).\\

\blem\label{lambda}${}$\\
Let $\rho,\rho'\in\R_+^N$ be such that $\La(\rho)>0$ and $\rho'_i\le\rho_i$ for any $i\in\{1,\ds,N\}$.\\
Then, $\La(\rho')>0$.
\elem\

\bpf${}$\\
Suppose by contradiction $\La(\rho')\le0$, that is $\La_{\mcal I,x}(\rho')\le0$ for some $\mcal I,x$.\\
This cannot occur for $\mcal I=\{i\}$ because it would mean $\rho'_i\ge\fr{8\pi(1+\a_i(x))}{a_{ii}}$, so the same inequality would for $\rho_i$, hence $\La(\rho)\le\La_{\mcal I,x}(\rho)\le0$.\\
Therefore, there must be some $\mcal I,x$ such that $\La_{\mcal I,x}(\rho')\le0$ and $\La_{\mcal I\sm\{i\},x}(\rho')>0$ for any $i\in\mcal I$; this implies
\beqa
\nonumber0&<&\La_{\mcal I\sm\{i\},x}(\rho')-\La_{\mcal I,x}(\rho')\\
\nonumber&=&2\sum_{j\in\mcal I}a_{ij}\rho'_i\rho'_j-a_{ii}{\rho'_i}^2-8\pi(1+\a_i(x))\rho'_i\\
\label{cond}&<&\rho'_i\l(2\sum_{j\in\mcal I}a_{ij}\rho'_j-8\pi(1+\a_i(x))\r).
\eeqa
It will be not restrictive to suppose, from now on, $\rho_1'\le\rho_1$ and $\rho_i'=\rho_i$ for any $i\ge2$, since the general case can be treated by exchanging the indices and iterating.\\
Assuming this, we must have $1\in\mcal I$, therefore we obtain:
\beqy
0&<&\La_{\mcal I,x}(\rho)-\La_{\mcal I,x}(\rho')\\
&=&8\pi(1+\a_1(x))(\rho_1-\rho_1')-a_{11}\l({\rho_1'}^2-\rho_1^2\r)-2\sum_{j\in\mcal I\sm\{1\}}a_{1j}(\rho_1'-\rho_1)\rho_j\\
&=&(\rho_1-\rho'_1)\l(8\pi(1+\a_1(x))-a_{11}(\rho_1'+\rho_1)-2\sum_{j\in\mcal I\sm\{1\}}a_{1j}\rho_j\r)\\
&<&(\rho_1-\rho'_1)\l(8\pi(1+\a_1(x))-2\sum_{j\in\mcal I}a_{1j}\rho'_j\r),
\eeqy
which is negative by \eqref{cond}. We found a contradiction.
\epf\

We will now show that if the parameter $\rho$ lies in the interior of $E$ then not only the functional is bounded from below but it is coercive in the space of zero-average functions. In particular, this fact allows to deduce the \qm{if} part in Corollary \ref{ex} from Theorem \ref{subcrit}.\\
On the other hand, if $\rho$ belongs to the boundary of $E$, then the scenario is quite different.\\

\blem\label{coe}${}$\\
Suppose $\rho\in\mri E$. Then, there exists a constant $C=C_\rho$ such that
$$J_\rho(u)\ge\fr{1}C\sum_{i=1}^N\int_\Si|\n u_i|^2\mrm dV_g-C.$$
Moreover, $J_\rho$ admits a minimizer which solves \eqref{liousy}.
\elem\

\bpf${}$\\
Choosing $\d\in\l(0,\fr{d(\rho,\pa E)}{\sqrt N|\rho|}\r)$ one has $(1+\d)\rho\in E$, so
\beqy
J_\rho(u)&=&\fr{\d}{2(1+\d)}\sum_{i,j=1}^Na^{ij}\int_\Si\n u_i\cd\n u_j\mrm dV_g+\fr{1}{1+\d}J_{(1+\d)\rho}(u)\\
&\ge&\fr{\d}{2\t(1+\d)}\sum_{i=1}^N\int_\Si|\n u_i|^2\mrm dV_g-C,
\eeqy
hence we get the former claim.\\
To get the latter, we notice that, due to invariance by translation, any minimizer can be supposed to be in $\ol H^1(\Si)^N$; therefore, we can restrict $J_\rho$ to this subspace. Here, the above inequality implies coercivity, and it is immediate to see that $J_\rho$ is also lower semi-continuous, hence the existence of minimizers follows from direct methods of calculus of variations.
\epf\

\blem\label{ncoe}${}$\\
Suppose $\rho\in\pa E$. Then, there exists a sequence $\l\{u^n\r\}_{n\in\N}\sub H^1(\Si)^N$ such that
$$\sum_{i=1}^N\int_\Si\l|\n u_i^n\r|^2\mrm dV_g\us{n\to+\infty}\lto+\infty\q\q\q\lim_{n\to+\infty}\fr{J_\rho\l(u^n\r)}{\sum_{i=1}^N\int_\Si\l|\n u_i^n\r|^2\mrm dV_g}\le0$$
\elem\

\bpf${}$\\
We first notice that $(1-\d)\rho\in E$ for any $\d\in(0,1)$. In fact, otherwise, from Lemma \ref{ee} we would get $\rho'\nin E$ as soon as $\rho'_i\ge(1-\d)\rho_i$ for some $i$, hence $\rho\nin\pa E$.\\
Now, suppose by contradiction that for any sequence $u^n$ one gets
$$\sum_{i=1}^N\int_\Si\l|\n u_i^n\r|^2\mrm dV_g\us{n\to+\infty}\lto+\infty\q\q\q\To\q\q\q\fr{J_\rho\l(u^n\r)}{\sum_{i=1}^N\int_\Si\l|\n u_i^n\r|^2\mrm dV_g}\ge\e>0.$$
Therefore, we would have
$$J_\rho(u)\ge\fr{\e}2\sum_{i=1}^N\int_\Si|\n u_i|^2\mrm dV_g-C;$$
hence, indicating as $\t'$ the smallest eigenvalue of the matrix $A$, for small $\d$ we would get
\beqy
J_\rho(u)&=&(1+\d)J_{(1+\d)\rho}(u)-\fr{\d}2\sum_{i,j=1}^Na^{ij}\int_\Si\n u_i\cd\n u_j\mrm dV_g\\
&\ge&\l((1+\d)\fr{\e}2-\fr{\d}{2\t'}\r)\sum_{i=1}^N\int_\Si|\n u_i|^2-C\\
&\ge&-C.\\
\eeqy
So we obtain $(1+\d)\rho\in E$; being also $(1-\d)\rho\in E$ (by Lemma \ref{ee}), we get a contradiction with $\rho\in\pa E$.
\epf\

To see what happens when $\rho\in\pa E$, we build an auxiliary functional using Lemma \ref{fanbn}.\\

\blem\label{jtilde}${}$\\
Fix $\rho'\in\pa E$ and define:
$$a_{\rho'}^n:=\fr{1}2\sum_{i,j=1}^Na^{ij}\int_\Si\n u_i^n\cd\n u_j^n\mrm dV_g\q\q\q b_{\rho'}^n:=J_{\rho'}\l(u^n\r)$$
$$J'_{\rho',\rho}(u)=J_\rho(u)-F_{\rho'}\l(\fr{1}2\sum_{i,j=1}^Na^{ij}\int_\Si\n u_i\cd\n u_j\mrm dV_g\r),$$
where $u^n$ is given by Lemma \ref{ncoe} and $F_{\rho'}$ by Lemma \ref{fanbn}.\\
If $\rho\in\mri E$, then $J'_{\rho',\rho}$ is bounded from below on $H^1(\Si)^N$ and its infimum is achieved by a solution of
$$-\D\l(u_i-\sum_{i,j=1}^Na^{ij}fu_j\r)=\sum_{j=1}^Na_{ij}\rho_j\l(\fr{\wt h_je^{u_j}}{\int_\Si\wt h_je^{u_j}\mrm dV_g}-1\r),\q\q\q i=1,\ds,N,$$
with $f=\l(F_{\rho'}\r)'\l(\fr{1}2\sum_{i,j=1}^Na^{ij}\int_\Si\n u_i\cd\n u_j\mrm dV_g\r)$.\\
On the other hand, $J'_{\rho',\rho'}$ is unbounded from below.
\elem\

\bpf${}$\\
For $\rho\in\mri E$, we can argue as in Lemma \ref{coe}, since the continuity follows from the regularity of $F$ and the coercivity from the behavior of $F'$ at the infinity.\\
For $\rho=\rho'$, if we take $u^n$ as in Lemma \ref{ncoe} we get
$$J'_{\rho',\rho'}\l(u^n\r)=b_{\rho'}^n-F_{\rho'}\l(a_{\rho'}^n\r)\us{n\to+\infty}\lto-\infty.$$
\epf\

Now we can prove the first half of Theorem \ref{subcrit}, that is $J_\rho$ is bounded from below if $\La(\rho)>0$.\\

\bpf[Proof of $\l\{\La>0\r\}\sub E$]${}$\\
Suppose by contradiction there is some $\rho'\in\pa E$ with $\La(\rho)>0$ and take a sequence $\rho^n\in E$ with $\rho^n\us{n\to+\infty}\lto\rho'$.\\
Then, by Lemma \ref{jtilde}, the auxiliary functional $J_{\rho',\rho^n}$ admits a minimizer $u^n$, so the functions $v_i^n$ defined as in \eqref{vi} solve
$$\l\{\bl-\D v_i^n=\sum_{j,j'=1}^Na_{ij}b^{jj',n}\l(\wt h_je^{v_j^n}-\rho_j^n\r)\\\int_\Si\wt h_i^ne^{v_i^n}\mrm dV_g=\rho_i^n\earr\r.\q\q\q i=1,\ds,N$$
where $b^{ij,n}$ is the inverse matrix of $b_{ij}^n:=\d_{ij}-a^{ij}f^n$, hence $b^{ij,n}\us{n\to+\infty}\lto\d_{ij}$.\\
We can then apply Theorem \ref{conc}. The first alternative is excluded, since otherwise we would get, for any $u\in H^1(\Si)^N$,
$$J'_{\rho',\rho'}(u)=\lim_{n\to+\infty}J'_{\rho',\rho^n}(u)\ge\lim_{n\to+\infty}J'_{\rho',\rho^n}\l(v^n\r)=J'_{\rho',\rho'}(v)>-\infty,$$
thus contradicting Lemma \ref{jtilde}.\\
Therefore, blow up must occur; this means, by Lemma \ref{linfty}, that $\s_i(p)\ne0$ for some $i\in\{1,\ds,N\}$ and some $p\in\Si$.\\
By Proposition \ref{lwzhang} follows $\La(\s)\le0$. On the other hand, since $\s_i\le\rho'_i$ for any $i$, Lemma \ref{lambda} yields $\La(\rho')\le0$, which contradicts our assumptions.
\epf\

To prove the unboundedness from below of $J_\rho$ in the case $\La(\rho)<0$ we will use suitable test functions, whose properties are described by the following:

\blem\label{test}${}$\\
Define, for $x\in\Si$ and $\la>0,\,\ph=\ph^{\la,x}$ as
$$\ph_i:=-2(1+\a_i(x))\log\max\{1,\la d(\cd,x)\}.$$
Then, as $\la\to+\infty$, one has
$$\int_\Si\n\ph_i\cd\n\ph_j\mrm dV_g=8\pi(1+\a_i(x))(1+\a_j(x))\log\la+O(1)$$
$$\ol{\ph_i}=-2(1+\a_i(x))\log\la+O(1)$$
$$\int_\Si\wt h_ie^{\sum_{j=1}^N\t_j\ph_j}\mrm dV_g\ge C\la^{-2(1+\a_i(x))}\q\q\q\tx{if}\q\sum_{i=1}^N\t_j(1+\a_j(x))>1+\a_i(x).$$
\elem\

\bpf${}$\\
It holds
$$\n\ph_i=\l\{\bll0&\tx{if }d(\cd,x)<\fr{1}\la\\-2(1+\a_i(x))\fr{\n d(\cd,x)}{d(\cd,x)}&\tx{if }d(\cd,x)>\fr{1}\la\earr\r..$$
Therefore, being $|\n d(\cd,x)|=1$ almost everywhere on $\Si$:
\beqy
&&\int_\Si\n\ph_i\cd\n\ph_j\mrm dV_g\\
&=&4(1+\a_i(x))(1+\a_j(x))\int_{\Si\sm B_\fr{1}\la(x)}\fr{\mrm dV_g}{d(\cd,x)^2}\\
&=&8\pi(1+\a_i(x))(1+\a_j(x))\log\la+O(1).
\eeqy\

For the average of $\ph_i$, we get
$$\int_\Si\ph_i\mrm dV_g=-2(1+\a_i(x))\int_{\Si\sm B_\fr{1}\la(x)}(\log\la+\log d(\cd,x))\mrm dV_g+O(1)=-2(1+\a_i(x))\log\la+O(1).$$\

For the last estimate, choose $r>0$ such that $\ol{B_\d(x)}$ does not contain any of the points $p_m$ for $m=1,\ds,M$, except possibly $x$.\\
Then, outside such a ball, $e^{\sum_{j=1}^N\t_j\ph_j}\le C\la^{-2\sum_{j=1}^N\t_j(1+\a_j(x))}$.\\
Therefore, under the assumptions of the Lemma,
$$\int_{\Si\sm B_\d(x)}\wt h_ie^{\sum_{i=1}^N\t_j\ph_j}\mrm dV_g=o\l(\la^{-2(1+\a_i(x))}\r),$$
hence
\beqy
&&\int_\Si\wt h_ie^{\sum_{i=1}^N\t_j\ph_j}\mrm dV_g\\
&\ge&\int_{B_\d(x)}\wt h_ie^{\sum_{i=1}^N\t_j\ph_j}\mrm dV_g\\
&\ge&C\l(\int_{B_\fr{1}\la(x)}d(\cd,x)^{2\a_i(x)}\mrm dV_g+\fr{1}{\la^{2\sum_{j=1}^N\t_j(1+\a_j(x))}}\int_{A_{\fr{1}\la,\d}(x)}d(\cd,x)^{2\a_i(x)-2\sum_{i=1}^N\t_j(1+\a_j(x))}\mrm dV_g\r)\\
&\ge&C\la^{-2(1+\a_i(x))},
\eeqy
which concludes the proof.
\epf\

\bpf[Proof of $E\sub\l\{\La\ge0\r\}$]${}$\\
Take $\rho,\mcal I,x$ such that $\La_{\mcal I,x}(\rho)<0$ and $\La_{\mcal I\sm\{i\},x}(\rho)\ge0$ for any $i\in\mcal I$, and consider the family of functions $\l\{u^\la\r\}_{\la>0}$ defined by
$$u_i^\la:=\sum_{j\in\mcal I}\fr{a_{ij}\rho_j}{4\pi(1+\a_i(x))}\ph_j^{\la,x}.$$
By Jensen's inequality we get
\beqy
J_\rho\l(u^\la\r)&\le&\fr{1}2\sum_{i,j=1}^Na^{ij}\int_\Si\n u_i^\la\cd\n u_j^\la\mrm dV_g+\sum_{i\in\mcal I}\rho_i\l(\ol{u_i^\la}-\log\int_\Si\wt h_ie^{u_i^\la}\mrm dV_g\r)+C\\
&=&\fr{1}2\sum_{i,j\in\mcal I}\fr{a_{ij}\rho_i\rho_j}{16\pi^2(1+\a_i(x))(1+\a_j(x))}\int_\Si\n\ph_i\cd\n\ph_j\mrm dV_g\\
&+&\sum_{i,j\in\mcal I}\fr{a_{ij}\rho_i\rho_j}{4\pi(1+\a_j(x))}\ol{\ph_j}-\sum_{i\in\mcal I}\rho_i\log\int_\Si\wt h_ie^{\sum_{j\in\mcal I}\fr{a_{ij}\rho_j}{4\pi(1+\a_j(x))}\ph_j}\mrm dV_g+C.
\eeqy
At this point, we would like to apply Lemma \ref{test} to estimate $J_\rho\l(u^\la\r)$. To be able to do this, we have to verify that 
$$\fr{1}{4\pi}\sum_{j\in\mcal I}a_{ij}\rho_j>1+\a_i(x)\q\q\q\fa\,i\in\mcal I.$$
If $\mcal I=\{i\}$, then $\rho_i>\fr{8\pi(1+\a_i(x))}{a_{ii}}$, so it follows immediately. For the other cases, it follows from \eqref{cond}.\\
So we can apply Lemma \ref{test} and we get from the previous estimates:
\beqy
J_\rho\l(u^\la\r)&\le&\l(\fr{1}{4\pi}\sum_{i,j\in\mcal I}a_{ij}\rho_i\rho_j-\fr{1}{2\pi}\sum_{i,j\in\mcal I}a_{ij}\rho_i\rho_j+2\sum_{i\in\mcal I}\rho_i(1+\a_i(x))\r)\log\la+C\\
&=&-\fr{\La_{\mcal I,x}(\rho)}{4\pi}\log\la+C\\
&\us{n\to+\infty}\lto&-\infty.
\eeqy
\epf\

\bpf[Proof of Corollary \ref{ex}]${}$\\
The coercivity in the case $\La<0$, hence the existence of minimizing solutions for \eqref{liousy} follows from Theorem \ref{subcrit} and Lemma \ref{coe}.\\
If instead $\La(\rho)\ge0$, then one can find out the lack of coercivity by arguing as before with the sequence $u^\la$, which verifies
$$\sum_{i=1}^N\int_\Si\l|\n u_i^\la\r|^2\mrm dV_g\us{\la\to+\infty}\lto+\infty\q\q\q J_\rho\l(u^\la\r)\le-\fr{\La_{\mcal I,x}(\rho)}{4\pi}\log\la+C\le C.$$
\epf\

\section{Proof of Theorem \ref{sharp}.}\

Here we will finally prove a sharp inequality in the case when the matrix $a_{ij}$ has non-positive entries outside its main diagonal.\\

As already pointed out in the introduction, the function $\La(\rho)$ can be written in a much shorter form under these assumptions, so the condition $\La(\rho)\ge0$ is equivalent to $\rho_i\le\fr{8\pi(1+\wt\a_i)}{a_{ii}}$ for any $i\in\{1,\ds,N\}$.\\
Moreover, thanks to Lemma \ref{ee}, in order to prove Theorem \ref{sharp} for all such $\rho$'s it will suffice to consider
\bequ\label{rho}
\rho^0:=\l(\fr{8\pi(1+\wt\a_1)}{a_{11}},\ds,\fr{8\pi(1+\wt\a_N)}{a_{NN}}\r).
\eequ
By what we proved in the previous Section, for any sequence $\rho^n\us{n\to+\infty}\nearrow\rho^0$ one has
$$\inf_{H^1(\Si)^N}J_{\rho^n}=J_{\rho^n}(u^n)\ge-C_{\rho^n},$$
so Theorem \ref{sharp} will follow by showing that, for a given sequence $\l\{\rho^n\r\}_{n\in\N}$, the constant $C_n=C_{\rho^n}$ can be chosen independently of $n$.\\

As a first thing, we provide a Lemma which shows the possible blow-up scenarios for such a sequence $u^n$.\\
Here, the assumption on $a_{ij}$ is crucial since it reduces largely the possible cases.\\

\blem\label{un}${}$\\
Let $\rho^0$ be as in \eqref{rho}, $\l\{\rho^n\r\}_{n\in\N}$ such that $\rho^n\nearrow\rho^0$, $u^n$ a minimizer of $J_{\rho^n}$ and $v^n$ as in \eqref{vi}. Then, up to subsequences, there exists a set $\mcal I\sub\{1,\ds,N\}$ such that:
\bite
\item If $i\in\mcal I$, then $\mcal S_i=\{x_i\}$ for some $x_i\in\Si$ which satisfy $\wt\a_i=\a_i(x_i)$ and $\s_i(x_i)=\rho_i^0$, and $v_i^n\us{n\to+\infty}\lto-\infty$ in $L^\infty_{\mrm{loc}}\l(\Si\sm\Cup_{j\in\mcal I}\{x_j\}\r)$.
\item If $i\nin\mcal I$, then $\mcal S_i=\es$ and $v_i^n\us{n\to+\infty}\lto v_i$ in $W^{2,q}_{\mrm{loc}}\l(\Si\sm\Cup_{j\in\mcal I}\{x_j\}\r)$ for some $q>1$ and some suitable $v_i$.
\eite
Moreover, if $a_{ij}<0$ then $x_i\ne x_j$.
\elem\

\bpf${}$\\
From Theorem \ref{conc} we get a $\mcal I\sub\{1,\ds,N\}$ such that $\mcal S_i\ne\es$ for $i\in\mcal I$.\\
If $\mcal S_i\ne\es$, then by Corollary \ref{equiv} one gets
$$0<\s_i(x)\le\rho_i^0\le\fr{8\pi(1+\a_i(x))}{a_{ii}}$$
for all $x\in\mcal S_i$, hence
\beqa\nonumber0&=&\La_{\{1,\ds,N\},x}(\s(x))\\
\label{lambdax}&\ge&\sum_{j=1}^N\l(8\pi(1+\a_j(x))\s_j(x)-a_{jj}\s_j(x)^2\r)\\
\nonumber&\ge&8\pi(1+\a_i(x))\s_i(x)-a_{ii}\s_i(x)^2\\
\nonumber&\ge&0.
\eeqa
Therefore, all these inequalities must actually be equalities.\\
From the last, we have $\s_i(x)=\rho_i^0=\fr{8\pi(1+\a_i(x))}{a_{ii}}$, hence $\a_i(x)=\wt\a_i$. On the other hand, since $\sum_{x\in\mcal S_i}\s_i(x)\le\rho_i^0$, it must be $\s_i(x)=0$ for all but one $x_i\in\mcal S_i$, so Corollary \ref{equiv} yields $\mcal S_i=\{x_i\}$.\\

Let us now show that $v_i^n\us{n\to+\infty}\lto-\infty$ in $L^\infty_{\mrm{loc}}$.\\
Otherwise, Theorem \ref{conc} would imply $v_i^n\us{n\to+\infty}\lto v_i$ almost everywhere, therefore by Fatou's Lemma we would get the following contradiction:
$$\s_i(x_i)<\int_\Si\wt h_ie^{v_i}\mrm dV_g+\s_i(x_i)\le\int_\Si\wt h_i^ne^{v_i^n}\mrm dV_g=\rho_i^n\le\rho_i=\s_i(x_i).$$
Since also inequality \eqref{lambdax} has to be an equality, we get $a_{ij}\s_i(x_i)\s_j(x_i)$ for any $i,j\in\mcal I$, so whenever $a_{ij}<0$ there must be $\s_j(x_i)=0$, so $x_i\ne x_j$.\\
Finally, if $\mcal S_i=\es$, the convergence in $W^{2,q}_{\mrm{loc}}$ follows from what we just proved and Theorem \ref{conc}.
\epf\

We basically showed that if a component of the sequence $v^n$ blows up, then all its mass concentrates at a single point which has the lowest singularity coefficient.\\
The next Lemma gives some more important information about the convergence or the blow-up of the components of $v^n$.\\

\blem\label{blowup}${}$\\
Let $v_i^n,\,v_i,\,\rho^0,\,\mcal I$ and $x_i$ as in Lemma \ref{un}.\\
Then,
\bite
\item If $i\in\mcal I$, then the sequence $v_i^n-\ol{v_i^n}$ converges to some $G_i$ in $W^{2,q}_{\mrm{loc}}\l(\Si\sm\Cup_{j\in\mcal I}\{x_j\}\r)$ for some $q>1$ and weakly in $W^{1,q'}(\Si)$ for any $q'\in(1,2)$, and $G_i$ solves:
$$\l\{\bl-\D G_i=\sum_{j\in\mcal I}a_{ij}\rho_j^0\l(\d_{x_j}-1\r)+\sum_{j\nin\mcal I}a_{ij}\l(\wt h_je^{v_j}-\rho_j^0\r)\\\ol{G_i}=0\earr\r..$$
\item If $i\nin\mcal I$, then $v_i^n\us{n\to+\infty}\lto v_i$ in the same space, and $v_i$ solves:
\bequ\label{eqv}
\l\{\bl-\D v_i=\sum_{j\in\mcal I}a_{ij}\rho_j^0\l(\d_{x_j}-1\r)+\sum_{j\nin\mcal I}a_{ij}\l(\wt h_je^{v_j}-\rho_j^0\r)\\\int_\Si\wt h_ie^{v_i}\mrm dV_g=\rho_i^0\earr\r..
\eequ
\eite
\elem\

\bpf${}$\\
From Lemma \ref{un} follows that, for $i\in\mcal I$, $\wt h_i^ne^{v_i^n}\us{n\to\infty}\wk\rho_i^0\d_{x_i}$ in the sense of measures; in fact, for any $\phi\in C(\Si)$
\beqy
\l|\int_\Si\wt h_i^ne^{v_i^n}\phi\mrm dV_g-\rho_i^0\phi(x_i)\r|&\le&\int_\Si\wt h_i^ne^{v_i^n}|\phi-\phi(x_i)|\mrm dV_g+\l|\rho_i^n-\rho_i^0\r||\phi(x_i)|\\
&\le&\e\int_{B_\d(x_i)}\wt h_i^ne^{v_i^n}\mrm dV_g+2\|\phi\|_{L^\infty(\Si)}\int_{\Si\sm B_\d(x_i)}\wt h_i^ne^{v_i^n}\mrm dV_g\\
&+&\l|\rho_i^n-\rho_i^0\r|\|\phi\|_{L^\infty(\Si)}\\
&\le&\e\rho_i^n+2\|\phi\|_{L^\infty(\Si)}o(1)+o(1)\|\phi\|_{L^\infty(\Si)},
\eeqy
which is, choosing properly $\e$, arbitrarily small. Therefore, $v_i$ solves \eqref{eqv}.\\
On the other hand, if $q'\in(1,2)$, then $\fr{q'}{q'-1}>2$, so any function $\phi\in W^{1,\fr{q'}{q'-1}}(\Si)$ is actually continuous, hence
\beqy
\l|\int_\Si\n\l(v_i^n-\ol{v_i^n}-G_i\r)\cd\n\phi\mrm dV_g\r|&=&\l|\int_\Si\l(-\D v_i^n+\D G_i\r)\phi dV_h\r|\\
&\le&\sum_{j\in\mcal I}a_{ij}\l|\int_\Si\wt h_je^{v_j^n}\phi\mrm dV_g-\rho_j^0\phi(p)\r|\\
&+&\sum_{j\nin\mcal I}a_{ij}\l|\int_\Si\wt h_j\l(e^{v_j^n}-e^{v_j}\r)\phi\mrm dV_g\r|\\
&\us{n\to+\infty}\lto&0.
\eeqy
Therefore, we get weak convergence in $W^{1,q'}(\Si)$ for any $q'\in(1,2)$; standard elliptic estimates yield convergence in $W^{2,q}_{\mrm{loc}}\l(\Si\sm\Cup_{j\in\mcal I}\{x_j\}\r)$.\\
In the same way we prove the same convergence of $v_i^n$ to $v_i$.
\epf\

From these information about the blow-up profile of $v^n$ we deduce an important fact which will be used to prove the main Theorem:\\

\bcor\label{wn}${}$\\
Let $v^n$ and $x_i$ be as in Lemmas \ref{un} and \ref{blowup} and $w^n$ be defined by $w_i^n=\sum_{j=1}^Na^{ij}v_j^n$ for $i\in\{1,\ds,N\}$.\\
Then, $w_i^n-\ol{w_i^n}$ is uniformly bounded in $W^{2,q}_{\mrm{loc}}(\Si\sm\{x_i\})$ for some $q>1$ if $i\in\mcal I$, whereas if $i\nin\mcal I$ it is bounded in $W^{2,q}(\Si)$.
\ecor\

\bpf${}$\\
Since $-\D w_i^n=\wt h_i^ne^{v_i^n}-\rho_i^n$, the claim follows from the boundedness of $e^{v_i^n}$ in $L^\infty_{\mrm{loc}}(\Si\sm\{x_i\})$ and from standard elliptic estimates.
\epf\

The last Lemma we need is a localized scalar Moser-Trudinger inequality for the blowing-up sequence.\\

\blem\label{mtloc}${}$\\
Let $w_i^n$ be as in Corollary \ref{wn} and $x_i$ as in the previous Lemmas. Then, for any $i\in\mcal I$ and any small $r>0$ one has
$$\fr{a_{ii}}2\int_{B_r(x_i)}\l|\n w_i^n\r|^2\mrm dV_g-\rho_i^n\l(\log\int_{B_r(x_i)}\wt h_ie^{a_{ii}w_i^n}\mrm dV_g-a_{ii}\ol{w_i^n}\r)\ge-C_r.$$
\elem\

\bpf${}$\\
Since $\Si$ is locally conformally flat, we can choose $r$ small enough so that we can apply Corollary \ref{as} up to modifying $\wt h_i^n$. We also take $r$ so small that $\ol{B_r(x_i)}$ contains neither any $x_j$ for $x_j\ne x_i$ nor any $p_m$ for $m=1,\ds,M$ (except possibly $x_i$).\\
Let $z^n$ be the solution of
$$\l\{\bll-\D z_i^n=\wt h_i^ne^{v_i^n}-\rho_i^n&\tx{in }B_r(x_i)\\z_i^n=0&\tx{on }\pa B_r(x_i)\earr\r..$$
Then, $w_i^n-\ol{w_i^n}-z_i^n$ is harmonic and it has the same value as $w_i^n-\ol{w_i^n}$ on $\pa B_r(x_i)$, so from standard estimates
$$\l\|w_i^n-\ol{w_i^n}-z_i^n\r\|_{C^1(B_r(x_i))}\le C\l\|w_i^n-\ol{w_i^n}\r\|_{C^1(\pa B_r(x_i))}\le C.$$
From Lemma \ref{blowup} we get
\beqy
\l|\int_{B_r(x_i)}\l|\n w_i^n\r|^2\mrm dV_g-\int_{B_r(x_i)}\l|\n z_i^n\r|^2\mrm dV_g\r|&=&\l|\int_{B_r(x_i)}\l|\n\l(w_i^n-z_i^n\r)\r|^2\mrm dV_g\r.\\
&+&\l.2\int_{B_r(x_i)}\n w_i^n\cd\n\l(w_i^n-z_i^n\r)\mrm dV_g\r|\\
&\le&\int_{B_r(x_i)}\l|\n\l(w_i^n-z_i^n\r)\r|^2\mrm dV_g\\
&+&2\l\|\n w_i^n\r\|_{L^1(\Si)}\l\|\n\l(w_i^n-z_i^n\r)\r\|_{L^\infty(B_r(x_i))}\\
&\le&C_r.
\eeqy
Moreover,
\beqy\int_{B_r(x_i)}\wt h_ie^{a_{ii}\l(w_i^n-\ol{w_i^n}\r)}\mrm dV_g&\le& e^{a_{ii}\l\|w_i^n-\ol{w_i^n}-z_i^n\r\|_{L^\infty\l(B_r\l(x_i\r)\r)}}\int_{B_r(x_i)}\wt h_ie^{a_{ii}z_i^n}\mrm dV_g\\
&\le&C_r\int_{B_r(x_i)}d(\cd,x_i)^{2\wt\a_i}e^{a_{ii}z_i^n}\mrm dV_g.
\eeqy
Therefore, since $\wt\a_i\le0$ and $a_{ii}\rho_i^n\le8\pi(1+\wt\a_i)$, we can apply Corollary \ref{as} to get the claim:
\beqy
\fr{a_{ii}}2\int_{B_r(x_i)}\l|\n w_i^n\r|^2\mrm dV_g-\rho_i^n\log\int_{B_r(x_i)}\wt h_ie^{a_{ii}\l(w_i^n-\ol{w_i^n}\r)}\mrm dV_g&\ge&\fr{1}{2a_{ii}}\int_{B_r(x_i)}\l|\n\l(a_{ii}z_i^n\r)\r|^2\mrm dV_g\\
&-&\rho_i^n\log\int_{B_r(x_i)}d(\cd,x_i)^{2\wt\a_i}e^{a_{ii}z_i^n}\mrm dV_g-C_r\\
&\ge&-C_r
\eeqy
\epf\

\bpf[Proof of Theorem \ref{sharp}]${}$\\
As noticed before, it suffices to prove the boundedness from below of $J_{\rho^n}\l(u^n\r)$ for a sequence $\rho^n\us{n\to+\infty}\nearrow\rho^0$ and a sequence of minimizers $u^n$ for $J_{\rho^n}$. Moreover, due to invariance by addition of constants, one can consider $v^n$ in place of $u^n$.\\
Let us start by estimating the term involving the gradients.\\
From Corollary \ref{wn} we deduce that the integral of $|\n w_i^n|^2$ outside a neighborhood of $x_i$ is uniformly bounded for any $i\in\mcal I$, and the integral on the whole $\Si$ is bounded if $i\nin\mcal I$.\\
For the same reason, the integral of $a_{ij}\n w_i^n\cd\n w_j^n$ on the whole surface is uniformly bounded. In fact, if $a_{ij}\ne0$, then $x_i\ne x_j$, then
\beqy
\l|\int_\Si\n w_i^n\cd\n w_j^n\mrm dV_g\r|&\le&\int_{\Si\sm B_r(x_j)}\l|\n w_i^n\cd\n w_j^n\r|\mrm dV_g+\int_{\Si\sm B_r(x_i)}\l|\n w_i^n\cd\n w_j^n\r|\mrm dV_g\\
&\le&\l\|\n w_i^n\r\|_{L^{q'}(\Si)}\l\|\n w_j^n\r\|_{L^{q''}\l(\Si\sm B_r\{x_j\}\r)}+\l\|\n w_i^n\r\|_{L^{q''}\l(\Si\sm B_r\{x_i\}\r)}\l\|\n w_j^n\r\|_{L^{q'}(\Si)}\\
&\le&C_r,
\eeqy
with $q$ as in Corollary \ref{wn}, $q'=\l\{\bll\fr{2q}{3q-2}<2&\tx{if }q<2\\1&\tx{if }q\ge2\earr\r.$ and $q''=\l\{\bll\fr{2q}{2-q}&\tx{if }q<2\\\infty&\tx{if }q\ge2\earr\r.$.\\
Therefore, we can write
\beqy
\sum_{i,j=1}^Na^{ij}\int_\Si\n v_i^n\cd\n v_j^n\mrm dV_g&=&\sum_{i,j=1}^Na_{ij}\int_\Si\n w_i^n\cd\n w_j^n\mrm dV_g\\
&\ge&\sum_{i\in\mcal I}a_{ii}\int_{B_r(x_i)}\l|\n w_i^n\r|^2\mrm dV_g-C_r.
\eeqy\

To deal with the other term in the functional, we use the boundedness of $w_i^n$ away from $x_i$: choosing $r$ as in Lemma \ref{mtloc}, we get
\beqy\int_\Si\wt h_i^ne^{v_i^n-\ol{v_i^n}}\mrm dV_g&\le&2\int_{B_r(x_i)}\wt h_i^ne^{v_i^n-\ol{v_i^n}}\mrm dV_g\\
&=&2\int_{B_r(x_i)}\wt h_ie^{\sum_{j=1}^Na_{ij}\l(w_j^n-\ol{w_j^n}\r)}\mrm dV_g\\
&\le&C_r\int_{B_r(x_i)}\wt h_ie^{a_{ii}\l(w_i^n-\ol{w_i^n}\r)}\mrm dV_g.
\eeqy
Therefore, using Lemma \ref{mtloc} we obtain
\beqy
J_{\rho^n}\l(v^n\r)&=&\fr{1}2\sum_{i,j=1}^Na^{ij}\int_\Si\n v_i^n\cd\n v_j^n\mrm dV_g-\sum_{i=1}^N\rho_i^n\l(\log\int_\Si\wt h_i^ne^{v_i^n}\mrm dV_g-\ol{v_i^n}\r)\\
&\ge&\sum_{i\in\mcal I}\l(\fr{a_{ii}}2\int_{B_r(x_i)}\l|\n w_i^n\r|^2\mrm dV_g-\rho_i^n\l(\log\int_{B_r(x_i)}\wt h_ie^{a_{ii}w_i^n}\mrm dV_g-a_{ii}\ol{w_i^n}\r)\r)-C_r\\
&\ge&-C_r
\eeqy

Since the choice of $r$ does not depend on $n$, the proof is complete.
\epf\

\section*{Acknowledgments}\

The author would like to express his gratitude Professor Andrea Malchiodi for the support and for the discussions about this topic.\\

\bibliography{finale}
\bibliographystyle{abbrv}

\end{document}